\magnification=1100
\baselineskip=12pt
\hsize12cm
\vsize18cm
\font\twelverm=cmr12
\font\twelvei=cmmi12
\font\twelvesy=cmsy10
\font\twelvebf=cmbx12
\font\twelvett=cmtt12
\font\twelveit=cmti12
\font\twelvesl=cmsl12

\font\ninerm=cmr9
\font\ninei=cmmi9
\font\ninesy=cmsy9
\font\ninebf=cmbx9
\font\ninett=cmtt9
\font\nineit=cmti9
\font\ninesl=cmsl9

\font\eightrm=cmr8
\font\eighti=cmmi8
\font\eightsy=cmsy8
\font\eightbf=cmbx8
\font\eighttt=cmtt8
\font\eightit=cmti8
\font\eightsl=cmsl8

\font\sixrm=cmr6
\font\sixi=cmmi6
\font\sixsy=cmsy6
\font\sixbf=cmbx6

\catcode`@=11 
\newskip\ttglue
\def\twelvepoint{\def\rm{\fam0\twelverm}
\textfont0=\twelverm  \scriptfont0=\ninerm  
\scriptscriptfont0=\sevenrm
\textfont1=\twelvei  \scriptfont1=\ninei  \scriptscriptfont1=\seveni
\textfont2=\twelvesy  \scriptfont2=\ninesy  
\scriptscriptfont2=\sevensy
\textfont3=\tenex  \scriptfont3=\tenex  \scriptscriptfont3=\tenex
\textfont\itfam=\twelveit  \def\it{\fam\itfam\twelveit}%
\textfont\slfam=\twelvesl  \def\sl{\fam\slfam\twelvesl}%
\textfont\ttfam=\twelvett  \def\tt{\fam\ttfam\twelvett}%
\textfont\bffam=\twelvebf  \scriptfont\bffam=\ninebf
\scriptscriptfont\bffam=\sevenbf  \def\bf{\fam\bffam\twelvebf}%
\tt  \ttglue=.5em plus.25em minus.15em
\normalbaselineskip=15pt
\setbox\strutbox=\hbox{\vrule height10pt depth5pt width0pt}%
\let\sc=\tenrm  \let\big=\twelvebig  \normalbaselines\rm}

\def\tenpoint{\def\rm{\fam0\tenrm}
\textfont0=\tenrm  \scriptfont0=\sevenrm  \scriptscriptfont0=\fiverm
\textfont1=\teni  \scriptfont1=\seveni  \scriptscriptfont1=\fivei
\textfont2=\tensy  \scriptfont2=\sevensy  \scriptscriptfont2=\fivesy
\textfont3=\tenex  \scriptfont3=\tenex  \scriptscriptfont3=\tenex
\textfont\itfam=\tenit  \def\it{\fam\itfam\tenit}%
\textfont\slfam=\tensl  \def\sl{\fam\slfam\tensl}%
\textfont\ttfam=\tentt  \def\tt{\fam\ttfam\tentt}%
\textfont\bffam=\tenbf  \scriptfont\bffam=\sevenbf
\scriptscriptfont\bffam=\fivebf  \def\bf{\fam\bffam\tenbf}%
\tt  \ttglue=.5em plus.25em minus.15em
\normalbaselineskip=12pt
\setbox\strutbox=\hbox{\vrule height8.5pt depth3.5pt width0pt}%
\let\sc=\eightrm  \let\big=\tenbig  \normalbaselines\rm}

\def\ninepoint{\def\rm{\fam0\ninerm}
\textfont0=\ninerm  \scriptfont0=\sixrm  \scriptscriptfont0=\fiverm
\textfont1=\ninei  \scriptfont1=\sixi  \scriptscriptfont1=\fivei
\textfont2=\ninesy  \scriptfont2=\sixsy  \scriptscriptfont2=\fivesy
\textfont3=\tenex  \scriptfont3=\tenex  \scriptscriptfont3=\tenex
\textfont\itfam=\nineit  \def\it{\fam\itfam\nineit}%
\textfont\slfam=\ninesl  \def\sl{\fam\slfam\ninesl}%
\textfont\ttfam=\ninett  \def\tt{\fam\ttfam\ninett}%
\textfont\bffam=\ninebf  \scriptfont\bffam=\sixbf
\scriptscriptfont\bffam=\fivebf  \def\bf{\fam\bffam\ninebf}%
\tt  \ttglue=.5em plus.25em minus.15em
\normalbaselineskip=11pt
\setbox\strutbox=\hbox{\vrule height8pt depth3pt width0pt}%
\let\sc=\sevenrm  \let\big=\ninebig  \normalbaselines\rm}

\def\eightpoint{\def\rm{\fam0\eightrm}
\textfont0=\eightrm  \scriptfont0=\sixrm  \scriptscriptfont0=\fiverm
\textfont1=\eighti  \scriptfont1=\sixi  \scriptscriptfont1=\fivei
\textfont2=\eightsy  \scriptfont2=\sixsy  \scriptscriptfont2=\fivesy
\textfont3=\tenex  \scriptfont3=\tenex  \scriptscriptfont3=\tenex
\textfont\itfam=\eightit  \def\it{\fam\itfam\eightit}%
\textfont\slfam=\eightsl  \def\sl{\fam\slfam\eightsl}%
\textfont\ttfam=\eighttt  \def\tt{\fam\ttfam\eighttt}%
\textfont\bffam=\eightbf  \scriptfont\bffam=\sixbf
\scriptscriptfont\bffam=\fivebf  \def\bf{\fam\bffam\eightbf}%
\tt  \ttglue=.5em plus.25em minus.15em
\normalbaselineskip=9pt
\setbox\strutbox=\hbox{\vrule height7pt depth2pt width0pt}%
\let\sc=\sixrm  \let\big=\eightbig  \normalbaselines\rm}

\def\twelvebig#1{{\hbox{$\textfont0=\twelverm\textfont2=\twelvesy
	\left#1\vbox to10pt{}\right.\n@space$}}}
\def\tenbig#1{{\hbox{$\left#1\vbox to8.5pt{}\right.\n@space$}}}
\def\ninebig#1{{\hbox{$\textfont0=\tenrm\textfont2=\tensy
	\left#1\vbox to7.25pt{}\right.\n@space$}}}
\def\eightbig#1{{\hbox{$\textfont0=\ninerm\textfont2=\ninesy
	\left#1\vbox to6.5pt{}\right.\n@space$}}}
 
\def\tit{\bigskip}

\def\Pun{\vbox{\hbox to 8.9pt {I\hskip-2.1pt P\hfil}}^1}
\def\Ptwo{\vbox{\hbox to 8.9pt {I\hskip-2.1pt P\hfil}}^2}
\def\Ptr{\vbox{\hbox to 8.9pt {I\hskip-2.1pt P\hfil}}^3}
\def\Pq{\vbox{\hbox to 8.9pt {I\hskip-2.1pt P\hfil}}^4}
\def\Pcq{\vbox{\hbox to 8.9pt {I\hskip-2.1pt P\hfil}}^5}
\def\Psx{\vbox{\hbox to 8.9pt {I\hskip-2.1pt P\hfil}}^6}
\def\Pn{\vbox{\hbox to 8.9pt {I\hskip-2.1pt P\hfil}}^n}
\def\Fn{\vbox{\hbox to 8.9pt {I\hskip-2.1pt F\hfil}}_n}
\def\epf{{\vrule height.9ex width.8ex depth-.1ex}}

\def\of{{\cal O}}
\def\i{{\cal I}}

\topinsert 
\endinsert
\font\big=cmbx10 scaled \magstep2

\centerline{\big \hbox {A note on rational surfaces in projective four-space.}}
\vskip 0.8truecm
\centerline {Ph. Ellia  \footnote {$^1$}
{Partially supported by MURST and Ferrara Univ. in the framework
of the project:
"Geometria algebrica, algebra commutativa e aspetti computazionali"}}
\medskip
{\eightpoint 
\centerline {Dipartimento di Matematica, Universit\`a di Ferrara}
\centerline {via Machiavelli 35 - 44100 Ferrara, Italy}
\centerline { e-mail: phe@dns.unife.it}}
\tit
\bigskip

\centerline {\bf 1. Introduction }
\tit
A few years ago Ellingsrud and Peskine proved ([12]) that there exist only finitely many irreducible components of the Hilbert scheme of $\Pq$ parametrizing smooth surfaces not of general type; in particular, as conjectured by Hartshorne and Lichtenbaum, the degree of smooth rational surfaces $S \subset \Pq$ is bounded. This result has been successively improved ([5], [8], [4], [9]) and today it is believed that if $S \subset \Pq$ is of non general type, then $deg(S) \leq 15$; also no rational surface of degree $d > 12$ is known.
\par
In this note we consider rational surfaces $S \subset \Pq$ ruled by cubics and quartics (i.e. possessing a base point free pencil of cubic or quartic rational curves) and we prove that such a surface has $deg(S) \leq 12$. (We recall that the classification of scrolls and conic bundles is known [3], [11], [6], [1]).
\par
The proof uses ad-hoc arguments which (unfortunately) do not seem to generalize.
\par
Using this result we then prove that if $S \subset \Pq$ is the image of a blow-up of $\Fn$ embedded by a linear system of the form $aC_0+bf-E_1-...-E_r$ (in the sequel, we will call such a linear system a "linear system on $\Fn$ with simple base points") then, again, $deg(S) \leq 12$.  
\tit
\centerline{\bf 2. Generalities }
\tit
Let $S \subset \Pq$ be a smooth, non-degenerated, rational surface. If $S$ is isomorphic to $\Ptwo$ then, by Severi's theorem, $S$ is a Veronese surface. If $S \simeq \Fn$ then $S$ is geometrically ruled and it is not difficult to see that $n=1$ and $S$ is a cubic scroll. Hence we may assume that $S$ is isomorphic to a blow-up of some $\Fn , n \geq 0$.
\tit
{\bf Definition 1.} {\it
We will say that $S$ is $a$-ruled if there exists on $S$ a base point free pencil of rational curves of degree $a$ in $\Pq$.}
\tit
{\bf Remark 1.} {\it
Such a pencil yields a morphism $p: S \to \Pun$ which presents $S$ as ruled by the curves of the pencil. Of course the same $S$ might be $a$-ruled for different values of $a$.
\par
Notice that since $S$ is not geometrically ruled, there is at least one singular fiber.}
\tit
{\bf Lemma 1.} {\it
Let $S \subset \Pq$ be a smooth, rational $a$-ruled surface, $a \geq 3$. If the general fiber of $p: S \to \Pun$ is degenerated in $\Pq$, then $S$ contains a plane curve of degree $d-a$, residual to a fiber in an hyperplane section.}
\tit 
{\it Proof:} Let $x$ be a general point of $\Pun$. The fiber $f_x$ is a smooth rational curve of degree $a$ in $\Pq$. By assumption $f_x$ is contained in an hyperplane, $H_x$ (note that $H_x$ is uniquely determined because $f_x$ is not a plane curve since $a \geq 3$). Let $C_x$ denote the residual curve: $C_x \sim H_x - f_x$. Since two general fibers are linearly equivalent, we have $C_x \sim C_y$ (they are both sections of $\of _S(1)\otimes p^*\of _{{\bf P}^1}(-1)$). Since $S$ is linearly normal (Severi's theorem) and since $f_x$ is not a plane curve, $h^0(\of _S(1-f_x))=1$. It follows that $C_x=C_y$. Now $C_x \subset H_x \cap H_y$, and since $S$ is non-degenerated, we may assume $H_x \neq H_y$, hence $C_x$ is a plane curve of degree $d-a$. \epf
\tit
The next proposition will be used several time in the sequel:
\tit
{\bf Proposition 2.} {\it
Let $S \subset \Pq$ be a smooth, non-degenerated, surface of degree $d$, not of general type. If $d \geq 9$, then $h^0(\i_S(3))=0$; in particular if $d > 9$ then $\pi \leq G(d,4)$ where $\pi$ is the sectional genus of $S$ and where $G(d,4)$ denotes the maximal genus of smooth degree $d$ curves in $\Ptr$ not lying on a cubic surface.}
\tit
{\it Proof:} See [10] \epf
\tit
{\bf Remark 2.} {\it
If $d > 12$, then $G(d,4) = 1+ {{d^2-3r(4-r)}\over {8}}$ where $d+r \equiv 0 (mod 4)$ and $0 \leq r <4$. In particular $\pi \leq 1+{{d^2}\over {8}}$; moreover if equality occurs then $\pi = G(d,4)$ and the general hyperplane section of $S$ is a.C.M. (arithmetically Cohen-Macaulay), but this is impossible because an a.C.M. surface in $\Pq$ not of general type has $d \leq 8$ (see [10]).
\par
In conclusion if $d >12$ and $S$ is not of general type then $\pi < 1+{{d^2}\over {8}}$.}
\tit
{\bf Corollary 3.} {\it
Let $S \subset \Pq$ be a smooth, $a$-ruled, rational surface. Assume $a \geq 3$. If the general fiber of $p: S \to \Pun$ is degenerated, then:
\par \noindent
(i) $\pi = {{(d-a-1)(d-a-2)}\over{2}}+a-1$.
\par \noindent
(ii) $1+2a-\sqrt{2a^2 -6a+5} \leq d \leq 1+2a+\sqrt{2a^2-6a+5}$.
\par \noindent
(iii) if $d > 12$, then ${{4a+6-2\sqrt{a^2-3a+15}}\over{3}} < d < {{4a+6+2\sqrt{a^2-3a+15}}\over{3}}$.}
\tit
{\it Proof:} (i) From lemma 1 it follows that $H \sim C+f$ where $C$ is a plane curve of degree $d-a$ and where $f$ is a rational curve of degree $a$. Since $a=f.H=f.C$, we get: $\pi =p_a(C \cup f)= p_a(C)+p_a(f)+a-1 = {{(d-a-1)(d-a-2)}\over{2}}+a-1$.
\par
(ii) The general hyperplane section of $S$ is non-degenerated in $\Ptr$ so its genus has to satisfy Castelnuovo's inequality: $\pi \leq ({{d}\over{2}}-1)^2$. Combining with (i) yields: $d^2+2d(-1-2a)+2a^2+10a-4 \leq 0$, and the result follows.
\par
(iii) By Remark 2: $\pi < 1+{{d^2}\over{8}}$, combining with (i) gives: $3d^2+2d(-4a-6)+4a^2+20a-8<0$, and we conclude. \epf
\tit
\centerline {\bf 3. $a$-ruled rational surfaces with $a \leq 3$.}
\tit
For sake of completeness we recall the following:
\tit
{\bf Proposition 4.} {\it
Let $S \subset \Pq$ be a smooth, non degenerated, rational surface.
\par
(i) if $S$ is a scroll ($a=1$), then $S$ is a cubic scroll.
\par
(ii) if $S$ is ruled in conics ($a=2$), then either $S$ is a Del Pezzo surface ($d=4$), or $S$ is a Castelnuovo surface ($d=5$).}
\tit
{\it Proof:} For (i) see [3], for (ii) see [11], [6] \epf
\tit
{\bf Proposition 5.} {\it
Let $S \subset \Pq$ be a smooth rational surface ruled in cubics ($a=3$).
\par
(i) $5 \leq d \leq 9$
\par
(ii) the possibilities for $(d, \pi)$ are: $(5,2),(6,3),(7,5),(8,8),(9,12)$.}
\tit
{\it Proof:} Since the fibers are cubics we can apply Corollary 3. From (ii) we get $5 \leq d \leq 9$, then we compute $\pi$ with (i). \epf
\tit
\centerline {\bf 4. Rational surfaces ruled in quartics.}
\tit
{\bf Lemma 6.} {\it
Let $S \subset \Pq$ be a smooth rational surface ruled in quartics. If the general fiber of $p: S \to \Pun$ is non-degenerated, then $h^1(\of_S(1))=0$ and $d \leq 9$.}
\tit
{\it Proof:} Consider Euler's sequence:
$$
0 \to M_S \to V\otimes \of _S \mathop\to ^{\rho} \of _S(1) \to 0
$$
($M:=\Omega_{{\bf P}^4}(1)$).
\par
We want to apply $p_*$ to this exact sequence. Restricting to a fiber we have:
$$
0 \to M_{f_x} \to V\otimes \of _{f_x} \mathop\to ^{\rho _x} \of _{f_x}(1) \to 0
$$
Notice that $h^0(\of _{f_x}(1))=5$ and $h^1(\of _{f_x}(1))=0$ for every $x$ in $\Pun$ (even if $f_x$ is singular); by base change it follows that $p_*(\of _S(1))$ is a rank 5 vector bundle on $\Pun$ and $R^ip_*(\of _S(1))=0, i>0$. Moreover, since for general $x$, $f_x$ spans $\Pq$, $\rho _x$ is an isomorphism and $h^0(M_{f_x})=0$ for general $x$. This implies $p_*(M_S)=0$ (it would be a torsion subsheaf of $p_*(V\otimes \of _S)=5.\of _{{\bf P}^1}$). Hence we get an injection: $0 \to 5.\of  _{{\bf P}^1} \to p_*(\of _S(1))$; let $T$ denote the cokernel, $T$ has finite support (it has rank zero). Taking cohomology in the exact sequence:
$$
0 \to 5.\of _{{\bf P}^1} \to p_*(\of _S(1)) \to T \to 0
$$
and since $h^0(p_*(\of _S(1))=h^0(\of _S(1))=5$ by Severi's theorem, we have $h^0(T)=0$, hence $T=0$ and $5.\of _{{\bf P}^1} \simeq p_*(\of _S(1))$. It follows that $h^1(p_*(\of _S(1))=0$. Since $R^ip_*(\of _S(1))=0,i>0$, by Leray's spectral sequence $h^1(\of _S(1))=h^1(p_*(\of _S(1))=0$ and $S$ is non-special.
\par
As shown in [2], non-special rational surfaces have $d \leq 9$. \epf
\tit
{\bf Remark 3.} {\it
Non-special rational surfaces are classified in [2].}
\tit
{\bf Proposition 7.} {\it
Let $S \subset \Pq$ be a smooth rational surface ruled in quartics, then $d \leq 12$.}
\tit
{\it Proof:} If the general fiber $f_x$ is a non-degenerated quartic in $\Pq$, we conclude with the previous proposition. If $f_x$ is degenerated, we conclude with Corollary 3. \epf
\tit
{\bf Remark 4.} {\it
As claimed in [7], every known rational surface contains a plane curve.}
\tit
\centerline {\bf Linear systems with simple base points on $F_n$.}
\tit
In this section we consider rational surfaces which are images of $\Fn$ by linear systems with simple base-points.
\medskip 
{\it Notations:} Let $S \subset \Pq$ be a smooth, non degenerated, surface isomorphic to $\Fn$ blown-up at $r$ points $y_1,...,y_r$.
\par
We have $Pic(\Fn )=C_0'{\bf Z}\oplus f'{\bf Z}$ where $(C_0')^2=-n$. Denoting by $C_0,f$ the strict transform of $C_0',f'$, we have $Pic(S)=C_0{\bf Z}\oplus f{\bf Z}\oplus E_1{\bf Z}\oplus ... \oplus E_r{\bf Z}$. We will work under the following assumptions:
\par
$$(*) \cases{ (a) & the $y_i$'s lie in different fibers of $\pi: \Fn \to \Pun$ \cr
(b) & If $n \geq 1$, no $y_i$ lies on $C_0'$ \cr
(c) & $H \sim aC_0+bf-E_1-...-E_r$ ("simple base points on $\Fn$")\cr} 
$$
\tit
{\bf Remark 5.} {\it
It follows that $S$ is $a$-ruled and that the fibers of the ruling $S \to \Pun$ have at most two irreducible components.}
\tit
The intersection theory on $S$ is given by: $C_0^2=-n,C_0E_i=0,C_0f=1,f^2=0,fE_i=0,E_iE_j=\delta_{ij}$.
\par \noindent
The canonical class is $K_S \sim -2C_0-(n+2)f+\Sigma E_i$.
\tit
We have the relations:
\par
1) $H^2=d$
\par
2) $2\pi -2=H(H+K)$
\par
3) $d(d-5)-10(\pi -1)+12\chi = 2K^2$
\par
\tit
After some computations we get:
\par
1) $d=-a^2n+2ab-r$
\par
2) $2\pi -2=-a^2n+an-2a+2ab-2b$
\par
3) $d(d-5)-10(\pi -1)=4-2r$
\tit
{\bf Lemma 8.} {\it
With notations as above, if $\pi < {{d^2}\over{8}}$, then $a \leq 9$.}
\tit
{\it Proof:} From 1): $r=-a^2n+2ab-d$, inserting in 3): $d^2-7d+3a^2n-5an+10a-4+b(10-6a)=0$, i.e.
$$
b={{d^2-7d+3a^2n-5an+10a-4}\over{6a-10}} \eqno (*)
$$
Using 2): $\pi -1=-{{an}\over{2}}(a-1)-a+{{(a-1)(d^2-7d+3a^2n-5an+10a-4)}\over{6a-10}}$
\par
Now, using this expression of $\pi -1$ in the inequality $\pi -1<{{d^2}\over {8}}$, yields $f_a(d) < 0$ (**), where:
$$
f_a(d)=d^2(a+1)-28(a-1)d+16a^2-16a+16
$$
Notice that $n$ has disappeared!
\par
We have ${\partial {f_a(d)}}\over {\partial d}$ = $0  \Leftrightarrow d={{14(a-1)}\over{a+1}}=:d_0$. Now $f_a(d_0)=(a-1)(16a - {196(a-1)\over {a+1}})+16$. If $a \geq 10$, we have $f_a(d) \geq f_a(d_0) > 0,\forall d$, contradicting (**). (indeed $(16a-{196(a-1)\over {a+1}} > 0$ if $a \geq 11$ and one checks directly that $f_{10}(d_0) > 0$.)
\par
In conclusion, if $\pi < 1+{d^2\over {8}}$ and if $a \geq 10$, then $f_a(d)>0,\forall d$, which contradicts (**) \epf
\tit
{\bf Lemma 9.} {\it
With notations as above, if $\pi < 1+{d^2\over {8}}$, then the possibilities are:
\par 
$a=5$: $d=11,6$
\par
$a=7$: $d=13,10$
\par
$a=8$: $d=7$
\par
or: $a \leq 4$.}
\tit
{\it Proof:} From lemma 8 we may assume $a \leq 9$ and the inequality $f_a(d) \leq 0$ (see proof of lemma 8); i.e. $d^2(a+1)-28(a-1)d+16a^2-16a+16 \leq 0$. Solving for the values of $a$ under consideration we obtain:
\par
$a=5,4 \leq d \leq 14$;
\par
$a=6, 5 \leq d \leq 15$;
\par
$a=7, 6 \leq d \leq 15$;
\par
$a=8, 7 \leq d \leq 15$;
\par
$a=9, 9 \leq d \leq 14$;
\par
On the other hand, using (*) of the proof of lemma 8:
\par
$$
(a-1)b={(a-1)(d^2-7d+10a-4)+(a-1)an(3a-5)\over {2(3a-5)}}
$$
$$
(a-1)b=n{a(a-1)\over {2}}+{(a-1)(d^2-7d+10a-4)\over {(6a-10)}}
$$
It follows that ${(a-1)(d^2-7d+10a-4)\over {(6a-10)}}$ is an integer. Now among the $(a,d)$ listed above, we take only those for which this further condition holds; this gives the statement of the lemma \epf 
\tit
{\bf Theorem 10.} {\it
Let $S \subset \Pq$ be a smooth, non degenerated, rational surface isomorphic to $\Fn$ blown-up at $r$ points $y_1,...,y_r$. Suppose assumptions (*) (see beginning of this section) are satisfied. Then $deg(S) \leq 12$.}
\tit
{\it Proof:} Assume $d>12$. By Remark 2, $\pi < 1+{d^2\over 8}$. By Lemma 9, $a \leq 4$ or $(a,d)=(7,13)$. In the first case, we know by Proposition 7 that $d \leq 12$. Let's consider the case $(a,d)=(7,13)$. We use relations 1), ...,3) before Lemma 8. From 2): $\pi -1=6b-7-21n$ (+); from 1): $-r =13+49n-14b$. Inserting in 3): $2b=7n+9$. Finally, from (+): $\pi = 21$. We observe that $21=G(13,4)$, hence arguing as in Remark 2, we conclude that $S$ is a.C.M.; but this is impossible ([10]) \epf
   
\tit
\centerline {\bf References}
\tit
\item{[1]} Abo H., Decker W., Sasakura N.: "An elliptic conic bundle in $\Pq$ arising from a stable rank three vector bundle", preprint (1997)
\item{[2]} Alexander, J.: "Surfaces non sp\'eciales dans $\Pq$", {\it Math. Zeitschrift}, {\bf 200}, 87-110 (1988) 
\item{[3]} Aure, A.: "On surfaces in projective fourspace", Thesis, Oslo (1987)
\item{[4]} Braun, R.-Cook, M.: "A smooth surface in $\Pq$ not of general type has degree at most 66", {\it Compositio Math.}, {\bf 107}, 1-9 (1997)
\item{[5]} Braun, R.-Fl\o ystad, G.: "A bound for the degree of smooth surfaces in $\Pq$ not of general type", {\it Compositio Math.}, {\bf 93}, 211-229 (1994)
\item{[6]} Braun, R.-Ranestad, Ch.: "Conic bundles in projective fourspace", in Algebraic Geometry: papers presented for the Europroj conferences in Catania and Barcelona, (Ed. P. Newstead)", {\it Lect. Notes in Pure and Applied Math.}, (M. Dekker Inc.) {\bf 200}, 331-339 (1998)
\item{[7]} Catanese, F.-Hulek, K.: "Rational surfaces in $\Pq$ containing a plane curve", {\it Ann. Mat. Pura ed Appl.}, (4) {\bf 172}, 229-256 (1997)
\item{[8]} Cook, M.: "An improved bound for the degree of smooth surfaces in $\Pq$ not of general type", {\it Compositio Math.}, {\bf 102}, 141-145 (1996)
\item{[9]} Cook, M.: "A smooth surface in $\Pq$ not of general type has degree at most 46", preprint.
\item{[10]} De Candia A.C.-Ellia, Ph.: "Some classes of non general type codimension two subvarieties in $\Pn$", {\it Ann. Univ. Ferrara}, vol XLIII, 135-156 (1997)
\item{[11]} Ellia, Ph.-Sacchiero, G.: "Smooth surfaces in $\Pq$ ruled in conics", in Algebraic Geometry: papers presented for the Europroj conferences in Catania and Barcelona, (Ed. P. Newstead)", {\it Lect. Notes in Pure and Applied Math.}, (M. Dekker Inc.) {\bf 200}, 49-62 (1998)
\item{[12]} Ellingsrud, G.-Peskine, Ch.: "Sur les surfaces lisses de $\Pq$", {\it Invent. Math.}, {\bf 95}, 1-11 (1989)

\end